\documentclass[12pt]{article}
\usepackage{graphicx}
\usepackage{graphicx}
\usepackage{amsmath, amsthm, amssymb, amscd}
\usepackage{setspace}
\usepackage{mathrsfs}

\newtheorem{lemma}{Lemma}
\newtheorem{theorem}{Theorem}
\newtheorem*{corollary 3.2 [4]}{Corollary 3.2 [7]}
\newtheorem*{corollary 1}{Corollary 1}
\newtheorem*{corollary 2}{Corollary 2}

\title{\textbf{Convex solutions to the power-of-mean curvature flow}\thanks{\textcopyright 2012 by
the author.}}

\author{Shibing Chen\thanks{Department of Mathematics, University of Toronto, Toronto, Ontario
Canada M5S 2E4 sbchen@math.toronto.edu.}}
\date{\today}

\begin{document}
 \maketitle

\begin{abstract}
 We prove some estimates for convex ancient solutions (the existence time for the solution starts from $-\infty$) to the power-of-mean curvature flow, when the power is strictly greater than $\frac{1}{2}$. As an application, we prove that in two dimension, the blow-down of the entire convex translating solution, namely $u_{h}=\frac{1}{h}u(h^{\frac{1}{1+\alpha}}x),$ locally uniformly converges to $\frac{1}{1+\alpha}|x|^{1+\alpha}$ as
 $h\rightarrow\infty$. Another application is that for generalized curve shortening flow (convex curve evolving in its normal direction with speed equal to a power of its curvature), if a convex compact ancient solution sweeps $\textbf{R}^{2}$, it has to be a shrinking circle. Otherwise the solution is defined in a strip region.
\end{abstract}

\section{\textbf{Introduction}}
Recently, classifying ancient convex solution to mean curvature flow has attracted much interest, due to its importance in studying the singularities
of mean curvature flow. Some important progress was made by Wang \cite{Wang}, and  Daskalopoulos, Hamilton and Sesum \cite{DHS}. In \cite{Wang} Wang proved that an entire convex translating solution to mean curvature flow must be rotationally symmetric which was a conjecture formulated explicitly by White in \cite{White}. Wang also constructed some entire convex translating solution with level set neither spherical
nor cylindrical in dimension greater or equal to 3. In the same paper, Wang also proved that if a convex ancient solution to the curve shortening flow sweeps the whole space $\textbf{R}^{2}$, it must be a shrinking circle, otherwise the convex ancient solution must be defined in a strip region and he indeed constructed such solutions by some compactness argument. Daskalopoulos, Hamilton and Sesum \cite{DHS} showed that besides the shrinking circle, the so called ¡±Angenent oval¡±(a convex ancient solution of the curve shortening flow discovered by Angenent that decomposes into two translating solutions of the flow) is the only other embedded convex compact ancient solution of the curve shortening flow. That means the corresponding curve shortening solution defined in a strip region constructed by Wang is exactly the ``Angenent oval".

The power-of-mean curvature flow, in which a hypersurface evolves in its normal direction with speed equal to a power $\alpha$ of its mean curvature $H$, was studied by Andrews \cite{An1}, \cite{An2}, \cite{An3}, Schulze \cite{Sc}, Chou and Zhu \cite{CZ} and  Sheng and Wu \cite{SW}
. Schulze \cite{Sc} called it $H^{\alpha}$-flow.
In the following, we will also call the one dimensional power-of-curvature flow the generalized curve shortening flow. Similar to the mean curvature flow, when one blows up the flow near the type II singularity appropriately, a convex translating solution  will arise, see \cite{SW} for details. It will be very interesting if one could classify the ancient convex solutions. In this paper, we will use the method developed by Wang \cite{Wang} to study the geometric asymptotic behavior of ancient convex solutions to $H^{\alpha}$-flow. The general equation for $H^{\alpha}$-flow is $\frac{\partial F}{\partial t}=-H^{\alpha}\vec{v}$, where $F:M\times [0, T)\rightarrow \textbf{R}^{n+1}$ is a time-dependent embedding of the evolving hypersurface, $\vec{v}$ is the unit normal vector to the hypersurface $F(M, t)$ in $\textbf{R}^{n+1}$ and $H$ is its mean curvature. If the evolving hypersurface can be represented as a graph of a function $u(x,t)$ over some domain in $\textbf{R}^{n}$, then we can project the evolution equation to the $(n+1)$th coordinate direction of $\textbf{R}^{n+1}$ and the equation becomes
$$u_{t}=\sqrt{1+|Du|^{2}}\left(\text{div}(\frac{Du}{\sqrt{1+|Du|^{2}}})\right)^{\alpha}.$$ Then a translating solution to the $H^{\alpha}$-flow will satisfy the equation $$\sqrt{1+|Du|^{2}}(\text{div}(\frac{Du}{\sqrt{1+|Du|^{2}}}))^{\alpha}=1,$$
which is equivalent to the following equation (3) when $\sigma=1$,
\begin{eqnarray}
L_{\sigma}(u)&=&(\sqrt{\sigma+|Du|^{2}})^\frac{1}{\alpha}\text{div}(\frac{Du}{\sqrt{\sigma+|Du|^{2}}})\\
&=&(\sigma+|Du|^{2})^{\frac{1}{2\alpha}-\frac{1}{2}}
\displaystyle{\sum_{i,j=1}^{n}}(\delta_{ij}-\frac{u_{i}u_{j}}{\sigma+|Du|^{2}})u_{ij}\\
&=&1,
\end{eqnarray}
where $\sigma\in [0,1]$, $\alpha\in (\frac{1}{2}, \infty]$ is a constant, $n=2$ is the dimension of $\textbf{R}^{2}$. If $u$ is a convex solution of (3), then $u+t$, as a function of $(x,t)\in \textbf{R}^{2}\times \textbf{R}$,
is a translating solution to the flow
\begin{eqnarray}
u_{t}=\sqrt{\sigma+|Du|^{2}}(\text{div}(\frac{Du}{\sqrt{\sigma+|Du|^{2}}}))^{\alpha}.
\end{eqnarray}
When $\sigma=1$, equation (4) is the non-parametric power-of-mean curvature flow. When $\sigma=0$, the level set $\{u=-t\}$, where $-\infty<t<-\inf u$, evolves by the power-of-mean curvature.

In the following we will assume $\sigma\in [0,1]$, $\alpha\in (\frac{1}{2}, \infty]$ and the dimension $n=2$, although some of the estimates do hold in high dimension. The main results of this paper are the following theorems.

\begin{theorem}
Suppose $u$ be an entire convex solution of (3). Let $u_{h}(x)=h^{-1}u(h^{\frac{1}{1+\alpha}}x).$ Then $u_{h}$ locally uniformly converges to $\frac{1}{1+\alpha}|x|^{1+\alpha}$,
as $h\rightarrow\infty.$
\end{theorem}
\begin{theorem}
Let $u_{\sigma}$ be an entire convex solution of (3). Then $u_{0}(x)=\frac{1}{1+\alpha}|x|^{1+\alpha}$ up to a a translation of the coordinate system.
When $\sigma\in(0,1]$, if $|D^{2}u(x)|=O(|x|^{3\alpha-1-2\alpha\beta})\ \text{as}\ |x|\rightarrow\infty,$ for some constant $\beta$ satisfying $\frac{1}{2\alpha}<\beta<\min\{1,\frac{1+\alpha}{2\alpha}\},$ then $u_{\sigma}$ is rotationally symmetric after a suitable translation of the coordinate system.
\end{theorem}
\begin{corollary 1}
A convex compact ancient solution to the generalized curve shortening flow which sweeps $\textbf{R}^{2}$ must be a shrinking circle.
\end{corollary 1}

\emph{Remark 1}. The condition $\alpha>\frac{1}{2}$ is necessary for our results. One can consider the translating solution $v(x)$ to (3) with $\sigma=1$ in one dimension. In fact when $\alpha\leq\frac{1}{2}$, the translating solution $v(x)$ is a convex function defined on the entire real line (\cite{CZ} page 28). Then one can construct a function $u(x,y)=v(x)-y$ defined on the entire plane, and $u$ will satisfy (3) with $\sigma=0$ and it is obviously not rotationally symmetric. We can also let $u(x, y)=v(x)$, which is an entire solution to (3) with $\sigma=1$ and it is not rotationally symmetric.
When the dimension $n$ is higher than two, similar examples can be given: we can take an entire rotationally
symmetric solution $v(x)$ to (1) with dimension $n\geq 2$ and $\sigma=1$, and then again let $u(x, y)=v(x)-y$, here $y$ is the $(n+1)$th coordinate for $R^{n+1}$. It is easy to see that $u$ will satisfy (3) with $n$ replaced by $n+1$ and $\sigma=0$, and the level set of $u$ is neither a sphere nor a cylinder.

Before embarking on the argument, we would like to point out that this elementary construction can be used to give a slight simplification of Wang's proof for Theorem 2.1 in \cite{Wang}( corresponding
to our Corollary 2 for $\alpha=1$). Let $v_{\sigma}$ be an entire convex solution to (3) in dimension $n$ with $\sigma\in(0,1]$. Then $u(x, y)=v_{\sigma}(x)-\sqrt{\sigma}y$ will be an entire convex solution to (3) in dimension $n+1$ with $\sigma=0$. Hence if one has proved the estimate in Corollary 2 for $\sigma=0$ in all dimensions, the estimates for $\sigma\in(0,1]$ follows immediately from the above construction.

The remainder of the paper is divided into three sections. The first contains the proof of Theorem 1 and the first part of Theorem 2. The second section establish Corollary 1 and the last section completes the proof of Theorem 2.

\section{\textbf{Proof of Theorem 1}}
For a given constant $h>0$, we denote
\begin{eqnarray*}
&\Gamma_{h}&=\{x\in R^{n}: u(x)=h\},\\
&\Omega_{h}&=\{x\in R^{n}: u(x)<h\},
\end{eqnarray*}
so that $\Gamma_{h}$ is the boundary of $\Omega_{h}.$
Denote $\kappa$ as the curvature of the level curve $\Gamma_{h}$. We have
\begin{eqnarray}
L_{\sigma}(u)&=&(\sigma+u_{\gamma}^{2})^{\frac{1}{2\alpha}-\frac{1}{2}}(\kappa u_{\gamma}+\frac{\sigma u_{\gamma\gamma}}{\sigma+u_{\gamma}^{2}})\\
&\geq& \kappa u_{\gamma}^{\frac{1}{\alpha}}=L_{0}(u),
\end{eqnarray}
where $\gamma$ is the unit outward normal to $\Omega_{h}$, and $u_{\gamma\gamma}=\gamma_{i}\gamma_{j}u_{ij}$.

\begin{lemma}

Suppose $u$ is a complete convex solution of (3). Suppose $u(0)=0$ and $\inf\{|x|: x\in \Gamma_{1}\}$
is achieved at $x_{0}=(0, -\delta)\in\Gamma_{1}$, for some $\delta>0$ very small. Let $D_{1}$ be the projection of $\Gamma_{1}$
on the axis $\{x_{2}=0\}$. Then the interval $(-R, R)$ is contained in $D_{1}$ with
\begin{eqnarray}
R\geq C_{1}(-\log\delta-C_{2})^{\frac{\alpha}{\alpha+1}},
\end{eqnarray}
where $C_{1}, C_{2} >0$ are independent of $\delta$.
\end{lemma}

\emph{Proof.}
We will prove the lemma when $\frac{1}{2}<\alpha\leq 1$ and indicate the small change needed for the case $\alpha>1$.
Suppose locally around $x_{0}$, $\Gamma_{1}$ is given by $z_{2}=f(z_{1})$. Then $f$ is a convex function satisfying $f(0)=-\delta$ and $f'(0)=0$. Take $a>0$
be a constant such that $f'(b)=1$. To prove (7) it is enough to prove
\begin{eqnarray}
a\geq C_{1}(-\log\delta-C_{2})^{\frac{\alpha}{\alpha+1}}.
\end{eqnarray}
For any $z=(z_{1}, z_{2})\in \Gamma_{1}$, where $z_{1}\in [0,a]$, let $\xi=\frac{z}{|z|}$,
from $\cite{Wang}$ we have
\begin{eqnarray}
u_{\gamma}(z)\geq \frac{\sqrt{1+f'^{2}}}{z_{1}f'-z_{2}},
\end{eqnarray}
. Since $L_{0}u\leq 1$, we have
\begin{eqnarray}
\frac{f''}{(1+f'^{2})^{\frac{3}{2}}}\frac{(1+f'^{2})^{\frac{1}{2\alpha}}}{(z_{1}f'-z_{2})^{\frac{1}{\alpha}}}\leq \kappa
u_{\gamma}^{\frac{1}{\alpha}}\leq 1.
\end{eqnarray}
 Hence
\begin{eqnarray}
f''(z_{1})&\leq&(1+f'^{2})^{\frac{3}{2}-\frac{1}{2\alpha}}(z_{1}f'-z_{2})^{\frac{1}{\alpha}}\\
&\leq&
10z_{1}^{\frac{1}{\alpha}}f'+10\delta
\end{eqnarray}
where $z_{2}=f(z_{1})$ and $f'(z_{1})\leq 1$ for $z_{1}\in (0,b)$. The inequality from (11) to (12) is trivial when $z_{2}\geq 0.$
When $z_{2}\leq 0$, since $|z_{2}|\leq \delta$, we have either $z_{1}f'\leq \delta$ or $z_{1}f'>\delta$, for the former
$(z_{1}f'-z_{2})^{\frac{1}{\alpha}}\leq (2\delta)^{\frac{1}{\alpha}}\leq 4\delta$, for the latter $(z_{1}f'-z_{2})^{\frac{1}{\alpha}}\leq (2z_{1}f')^{\frac{1}{\alpha}}\leq 4z_{1}^{\frac{1}{\alpha}}f'$, since $f'(z_{1})\leq 1$.
We consider the equation
\begin{eqnarray}
\rho''(t)=
10t^{\frac{1}{\alpha}}\rho'+10\delta
\end{eqnarray}
with initial conditions $\rho(0)=-\delta$ and $\rho'(0)=0$. Then for $t\in(0,b)$ we have
\begin{eqnarray}
\rho'(t)=10\delta e^{\frac{10\alpha}{\alpha+1}t^{\frac{\alpha+1}{\alpha}}}\int_{0}^{t} e^{-\frac{10\alpha}{\alpha+1}s^{\frac{\alpha+1}{\alpha}}}ds.
\end{eqnarray}
Since $\int_{0}^{\infty} e^{-\frac{10\alpha}{\alpha+1}s^{\frac{\alpha+1}{\alpha}}}ds$ is bounded above by some constant $C$, we have
\begin{eqnarray*}
1=\rho'(a)&=&10\delta e^{\frac{10\alpha}{\alpha+1}a^{\frac{\alpha+1}{\alpha}}}\int_{0}^{a} e^{-\frac{10\alpha}{\alpha+1}s^{\frac{\alpha+1}{\alpha}}}ds.\\
&\leq& C_{1}\delta e^{\frac{10\alpha}{\alpha+1}a^{\frac{\alpha+1}{\alpha}}},
\end{eqnarray*}
from where (8) follows. When $\alpha>1$, we need only to introduce a number $c$ such that $f'(c)=\frac{1}{2}$ and then we can find the lower bound
of $a-c$ in a similar way.

\emph{Remark 2}. It follows from Lemma 1 that when $\delta$ is sufficiently small, by convexity and in view of Figure 1, we see that $\Omega_{1}$ contains the shadowed region. Then it is easy to check that
$\Omega_{1}$ contains an ellipse
\begin{eqnarray}
E=\{(x_{1}, x_{2})|\frac{x_{1}^{2}}{(\frac{R}{6})^{2}}+\frac{(x_{2}-\frac{7\delta^{*}-5\delta}{12})^{2}}{(\frac{\delta^{*}+\delta}{4})^{2}}=1\},
\end{eqnarray}
where $\delta^{*}$ is a positive constant such that $u(0,\delta^{*})=1$ and $R$ is defined in the Lemma 1.

\begin{figure}
  \includegraphics[width=0.8\linewidth]{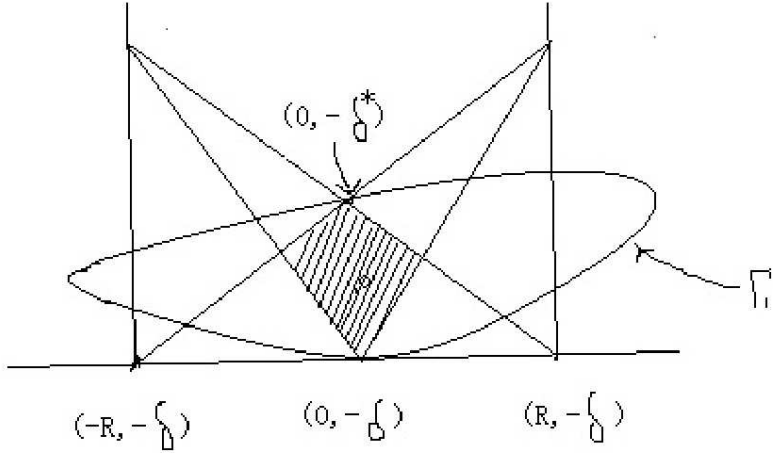}\\
  \caption{$\Gamma_{1}$ contains the shadow part.}\label{dsds}
\end{figure}

\emph{Remark 3.}
One can also establish similar lemma in higher dimensions, which says $D_{1}$ (convex set with dimension greater than 1) contains a ball centered at
the origin with radius $R\geq C_{n}(-\log\delta-C)^{\frac{\alpha}{\alpha+1}}$. For the details of how to reduce the situation to lower dimensional case we refer the reader to the proof of Lemma 2.6 in \cite{Wang}.

\begin{lemma}
Suppose $u$ is a complete convex solution of (3). Suppose $u(0)=0$, $\delta$ and $\delta^{*}$
are defined as in Lemma 1 and Remark 2. Then if $\delta$ and $\delta^{*}$  are small enough, $u$ is defined in a strip region.
\end{lemma}
The proof of Lemma 2 is based on a careful study of the shape of the level curve of $u$, we will give an important corollary first.

\begin{corollary 2}
Suppose $u$ is an entire convex solution of (3) in $\textbf{R}^{2}$, then
\begin{eqnarray}
u(x)\leq C(1+|x|^{1+\alpha}),
\end{eqnarray}
for some constant $C$ depending only on the upper bound of $u(0)$ and $|Du(0)|$.
\end{corollary 2}
\emph{Proof.}
By subtracting a constant we may assume $u(0)=0$. It is enough to prove that
$\text{dist}(0,\Gamma_{h})\geq C h^{\frac{1}{1+\alpha}}$ for all large $h$. By the rescaling $u_{h}(x)=\frac{1}{h}u(h^{\frac{1}{1+\alpha}}x)$
we need only to prove $\text{dist}(0,\Gamma_{1,u_{h}})\geq C$. Notice that $|Du_{h}(0)|=\frac{1}{h^{\frac{\alpha}{1+\alpha}}}|Du(0)|\rightarrow
0, \text{as}\ h\rightarrow \infty.$ Hence by convexity $\inf_{B_{r}(0)}u_{h}$ goes to 0 uniformly for fixed radius $r$. Note also that $u_{h}$
satisfies equation (3) with $\sigma \rightarrow 0\ \text{as}\ \ h\rightarrow \infty$.

If the estimate $$\text{dist}(0,\Gamma_{1,u_{h}})\geq C,\ \text{for all large}\ h$$ fails, we can find a
sequence $h_{k}\rightarrow\infty$ such that $\delta_{k}=\inf\{|x|: x\in\Gamma_{1,u_{h_{k}}}\}\rightarrow 0$.
Now, we take $\delta_{k}^{*}$ as in Remark 2 with respect to $u_{h_{k}}$.  $\delta_{k}^{*}$ has a positive lower bound $\delta^{*}$, otherwise
by Lemma 2 $u_{h_{k}}$ can not be an entire solution for large $k$.

If $\delta_{k}^{*}\leq 1000$ for all large $k$, since the ellipse $E_{k}$ defined for $u_{h_{k}}$ as in Remark 2 is contained in $\Omega_{1,u_{h_{k}}}$ and the distance between the center $O_{k}$ of $E_{k}$ and the origin
is bounded above by 1000, by the previous discussion we know $u_{h_{k}}(O_{k})$ is bounded bellow by $-1$ when $k$ is large. Let
$E_{k}(t)$ be the solution to the generalized curve shortening flow starting from time $t=-1$ , with initial data $E_{k}(-1)=E_{k}$.
(1) When $\sigma=0$, $\partial \Omega_{-t,u_{h_{k}}}$ evolves under the generalized curve shortening flow, we have
the inclusion $E_{k}(t)\subset\partial \Omega_{-t,u_{h_{k}}}$ for all $t>-1$. Hence $\inf_{B_{1000}(0)}u_{h_{k}}$ is smaller
than 1 minus the time needed for $E_{k}$ to shrink to $O_{k}$. However, by the size of  $E_{k}$, the time needed for it to shrink to a point goes to infinity as $k$ goes to infinity, which is contradictory to the discussion at the beginning of the proof that $u_{h_{k}}$ converges to 0 uniformly in the ball $B_{1000}(0)$ as $h_{k}$ goes to infinity. (2)When $\sigma\in(0, 1]$, we can take $v_{k}$ as the solution of $L_{\sigma_{k}}v=1$ in $E_{k}$ with $v=1$ on $\partial E_{k}$, where
$\sigma_{k}=h_{k}^{-\frac{2\alpha}{1+\alpha}}$. Passing to a subsequence and adjusting the size of $E_{k}$ if necessary, we can assume $E_{k}$
converge to some ellipse $E$ with the length of its long axis very large, the length of its short axis bigger than some fixed positive number and
the distance from its center to the origin is less than 1000.
Then $v_{k}$ converges to a solution of the generalized curve shortening flow, and a contradiction can be made as for the case $\sigma=0$

Otherwise, by the definition of $b$ in the proof of Lemma 1 and the convexity of $\Omega_{1,u_{h}}$ we
can find a disc $B_{k}$ with center $O=(0,50)$ and radius 20 inside $\Omega_{1,u_{h_{k}}}$, obviously it will take time more than 2 for $B_{k}$ to
shrink to $O$. We can take $B_{k}(t)$ as a solution to the generalized curve shortening flow starting from time $t=-1$
with $B_{k}(-1)=B_{k}$, then a similar contradiction
will be made as before.

\emph{Remark 4}. The estimate in Corollary 1 is also true for higher dimensions, one can prove it by reducing the problem to two dimensional case similar to the corresponding part in \cite{Wang}.

\emph{Proof of Lemma 2.} By rotating coordinates we assume the axial directions of $E$ in Remark 2 are the same with those of the coordinate system. Denote $M_{u}$ as the graph of $u$, and as in $\cite{Wang}$ we divide it into two parts, $M^{+}$ and $M^{-}$, where $M^{+}=\{(x,u(x))\in \textbf{R}^{3}: u_{x_{2}}\geq 0\}$ and $M^{-}=\{(x,u(x))\in \textbf{R}^{3}: u_{x_{2}}\leq 0\}$. Then $M^{\pm}$ are the graphs of functions $f^{\pm}$, namely the graph of $x_{2}=f^{\pm}(x_{1},x_{3}), (x_{1},x_{2})\in D$ and $D$ is the projection of $M_{u}$ on the $x_{1}x_{3}$ plane.
 The function $f^{+}$ is concave and $f^{-}$ is convex, and we have $x_{3}=u(x_{1}, f^{\pm}(x_{1}, x_{3}))$.
Let  \begin{eqnarray}
f=f^{+}-f^{-}.
\end{eqnarray}
Now it is easy to see that $f$ is positive and concave in $D$. Also note that $f$ is vanishing on $\partial{D}$. For any $h>0$ we also denote $f_{h}(x_{1})=f(x_{1},h),$
$f^{\pm}_{h}(x_{1})=f(x_{1}, h)$, and $D_{h}=\{x_{1}\in \textbf{R}^{1}: (x_{1}, h)\in D\}$. Then $f_{h}$ is a positive, concave function in $D_{h}$,
vanishing on $\partial{D_{h}}$, and $D_{h}=(-\underline{a}_{h}, \overline{a}_{h})$ is an interval containing the origin. We denote $b_{h}=f_{h}(0)$.
We will consider the case $\sigma=0$ first.

\emph{Claim 1:} suppose $h$ large, $f_{1}(0)=\delta^{*}+\delta$ small, $b_{h}\leq 4$ and $\underline{a}_{h}, \overline{a}_{h}\geq b_{h}$. Then
$\overline{a}_{h}\geq \frac{1}{1000}\frac{h}{b_{h}^{\alpha}}$ for $\alpha\leq 1$ and $\overline{a}_{h}\geq \frac{1}{1000}\frac{h^{\frac{1}{2\alpha-1}}}{b_{h}^{\frac{1}{2\alpha-1}}}$ for $\alpha>1$.

\emph{Proof.} Without loss of generality, we assume $\overline{a}_{h}\leq\underline{a}_{h}$. Denote $U_{h}=\Omega_{h}\cap\{x_{1}>0\}$.
Similar to that in $\cite{Wang}$, we have that the arc-length of the image of $\Gamma_{s}\cap\{x_{1}>0\}$ under Gauss map is bigger then $\frac{\pi}{6}.$ Notice that
$\Omega_{1}$ contains $E$, which was defined in Remark 2. When $\delta$ and $\delta^{*}$ are very small, $E$ is very thin and long. The centre of
$E$ is very close to the origin, in fact for our purpose we can just pretend $E$ is centered at the origin. By convexity of $\Omega_{h}$
and in view of Figure 2, we see that $\Gamma_{s}\cap\{x_{1}>0\}$ is trapped between two lines $\ell_{1}$ and $\ell_{2}$, and the slopes of  $\ell_{1}$ and $\ell_{2}$ are very close to 0 when $E$ is very long and thin. Then it is clear that the largest distance from the points on $\Gamma_{s}\cap\{x_{1}>0\}$
to the origin can not be bigger than $10\overline{a}_{h}$. By convexity of $u$, we have $u_{\gamma}(x)\geq \frac{h}{20\overline{a}_{h}},$
for $x\in\Gamma_{s}\cap\{x_{1}>0\}$. Since $\Gamma_{s}\cap\{x_{1}>0\}$ evolves under the generalized curve shortening flow, when $\alpha\leq 1$ we have the following estimate
\begin{eqnarray}
\frac{d}{ds}(|U_{s}|)&=&\int_{\Gamma_{s}\cap\{x_{1}>0\}}\kappa^{\alpha}d\xi\\
&=&\int_{\Gamma_{s}\cap\{x_{1}>0\}}u_{\gamma}^{\frac{1}{\alpha}-1}\kappa d\xi\\
&\geq& \frac{1}{50}(\frac{h}{\overline{a}_{h}})^{\frac{1}{\alpha}-1}\frac{\pi}{6},
\end{eqnarray}
from (18) to (19) we used the equation $\kappa u_{\gamma}^{\frac{1}{\alpha}}=1$. The claim follows by the simple fact
$\frac{3}{2}b_{h}\overline{a}_{h}\geq |U_{h}|\geq \frac{1}{50}(\frac{h}{\overline{a}_{h}})^{\frac{1}{\alpha}-1}\frac{\pi}{6}\frac{h}{2}$.

\begin{figure}
  \includegraphics[width=0.8\linewidth]{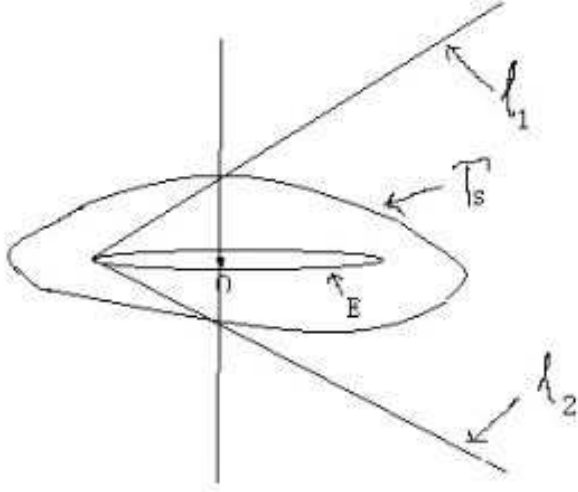}\\
  \caption{$\Gamma_{s}\cap\{x_{1}>0\}$ is trapped between two lines}\label{dsds}
\end{figure}

When $\alpha>1$, denote $l_{s}$ as the arc length of $\Gamma_{s}\cap\{x_{1}>0\}$, by the above discussion, it is not hard to see that
$l_{s}\approx C\overline{a}_{h}$. Then by a simple application of Jensen's inequality, we have
\begin{eqnarray*}
\frac{d}{ds}(|U_{s}|)&=&\int_{\Gamma_{s}\cap\{x_{1}>0\}}\kappa^{\alpha}d\xi\\
&=&\l_{s}\int_{\Gamma_{s}\cap\{x_{1}>0\}}\kappa^{\alpha}\frac{1}{l_{s}}d\xi\\
&\geq&l_{s}(\int_{\Gamma_{s}\cap\{x_{1}>0\}}\frac{\kappa}{l_{s}}d\xi)^{\alpha}\geq C l_{s}^{1-\alpha}\geq C\overline{a}_{h}^{1-\alpha},
\end{eqnarray*}
then by the simple fact that $\frac{3}{2}b_{h}\overline{a}_{h}\geq |U_{h}|$ we can finish the proof in the same way as the previous case.

\emph{
Claim 2:
}Denote $f_{k}=f_{h_{k}}$. Then
\begin{eqnarray}
f_{k}(0)\leq f_{k-1}(0)+C_{0}2^{\frac{-k}{C}} \ \text{for all}\  k \ \text{large},
\end{eqnarray}
where $C_{0}$ is a fixed constant, and $C$ depends only on $\alpha$.

Lemma 2 follows from Claim 1 and Claim 2 in the following way. Let the convex set $P$ be the projection of the graph of $g$ on the plane $\{x_{3}=0\}$, by Claim 2 and the fact that $P$ contains $x_{1}$-axis (it follows from Claim 1), $P$ must equal to $I\times \textbf{R}$ for some interval $I\subset [0,\displaystyle{\lim_{k\rightarrow\infty}} g_{k}(0)].$ Then, by (17) $\mathcal{M}_{u}$ is also contained in a strip region as stated in Lemma 2.

The proof of Claim 2 can be done by following the lines in $\cite{Wang}$ closely, but one should be very careful for choosing the proper constants in the proof which is very different from the case in $\cite{Wang}$.

\emph{Proof of Theorem 1 and the first part of Theorem 2.} First we prove that one can find a subsequence of $u_{h},$ where $u_{h}(x)=h^{-1}u(h^{\frac{1}{1+\alpha}}x),$ which converges to
$\frac{1}{1+\alpha}|x|^{1+\alpha}.$

By subtracting a constant we may suppose $u(0)=0.$ Suppose $x_{n+1}=b\cdot x$ is the tangent plane of $u$ at $0.$ By Corollary 2 and the convexity of $u$
we have $$b\cdot x\leq u(x)\leq C(1+|x|^{1+\alpha}).$$ Hence, $$h^{-\frac{\alpha}{1+\alpha}}b\cdot x\leq u_{h}(x)\leq C(\frac{1}{h}+|x|^{1+\alpha}).$$
It is easy to see that $Du_{h}$ is locally uniformly bounded. Hence $u_{h}$ sub-converges to a convex function $u_{0}$ which satisfies $u_{0}(0)=0,$
and $$0\leq u_{0}(x)\leq C|x|^{1+\alpha}.$$ Then it is easy to check that $u_{0}$ is an entire convex viscosity solution to equation (3) with $\sigma=0,$ and the comparison principle holds on any bounded domain. By using comparison principle it is easy to prove $\{x|u_{0}(x)=0\}=\{0\}.$

Now since $\{x|u_{0}(x)=0\}=\{0\},$  $\Gamma_{1,u_{0}}=\{x|u_{0}(x)=1\}$ is a bounded convex curve, and
the level set $\{x|u_{0}(x)=-t\}$, with time $t\in(-\infty, 0)$, evolves under the generalized curve shortening flow, from [1], [2] we have the following asymptotic behavior of the convex solution $u_{0}$ of $L_{0}u=1$
\begin{eqnarray}
 u_{0}(x)=\frac{1}{1+\alpha}|x|^{1+\alpha}+\varphi(x),\ \text{where}\ \varphi(x)=o(x^{1+\alpha}),\ \text{for}\ x\neq 0\ \text{near the origin}.
\end{eqnarray}
In fact, if the initial level curve is in a sufficiently small neighborhood of circle, by Lemma 4 in the beginning of the fourth section, we have that $|\varphi(x)|\leq C|x|^{1+\alpha+\eta}$ for some small positive $\eta,$ where $C$ is a constant depending only on the initial closeness to the circle. Hence, given any $\epsilon >0$, for small enough $h > 0$, we have $$B_{(1-\epsilon)r}(0)\subset \Omega_{h,u_{0}}\subset B_{(1+\epsilon)r}(0),$$ where $r=\left((1+\alpha)h\right)^{\frac{1}{1+\alpha}}.$ So there is a sequence $h_{m}\rightarrow\infty$ such that
$$B_{(1-\frac{1}{m})r_{m,i}}(0)\subset \Omega_{h_{m},u}\subset B_{(1+\frac{1}{m})r_{m,i}}(0),$$ where $r_{m,i}=\left((1+\alpha)ih_{m}\right)^{\frac{1}{1+\alpha}}, i=1,\cdots,m.$
Then $u_{h_{m}}$ sub-converges to $\frac{1}{1+\alpha}|x|^{1+\alpha}.$

Since $u_{0}$ is an entire convex solution to $L_{0}u=1,$ from the above argument, we can find a sequence $h_{m},$
such that $u_{0h_{m}}(x)=\frac{1}{h_{m}}u_{0}(h_{m}^{\frac{1}{1+\alpha}}x)$ locally uniformly converges to $\frac{1}{1+\alpha}|x|^{1+\alpha}.$ Hence, the
sublevel set $\Omega_{\frac{1}{1+\alpha},u_{0h_{m}}}$  satisfies $$B_{1-\epsilon_{m}}(0)\subset\Omega_{\frac{1}{1+\alpha},u_{0h_{m}}}\subset B_{1+\epsilon_{m}}(0),$$ where $\epsilon_{m}\rightarrow 0$ as $m\rightarrow\infty.$ By the discussion below (54), we have
$$u_{0h_{m}}(x)=\frac{1}{1+\alpha}|x|^{1+\alpha}+\varphi(x),$$ where $|\varphi(x)|\leq C|x|^{1+\alpha+\eta}$ for some fixed small positive $\eta,$ and the constant
$C$ is independent of $m.$ Replacing $x$ by $h_{m}^{-\frac{1}{1+\alpha}}x$ in the above asymptotic formula, we have
$$u_{0}(x)=\frac{1}{1+\alpha}|x|^{1+\alpha}+h_{m}\varphi(h_{m}^{-\frac{1}{1+\alpha}}x),$$ where for a given $x,$ $h_{m}\varphi(h_{m}^{-\frac{1}{1+\alpha}}x)\rightarrow 0.$
Hence $u_{0}(x)=\frac{1}{1+\alpha}|x|^{1+\alpha}.$  So we have proved Theorem 1 and the first part of Theorem 2.

\section{\textbf{Proof of Corollary 1}}
We will follow the lines in the section 4 of \cite{Wang}. It will be accomplished by the following lemma which is also true for higher dimensions, but we will only state it for $\textbf{R}^{2}$.
\begin{lemma}
Suppose $\Omega$ is a smooth, convex, bounded domain in $\textbf{R}^{2}$. Let $u$ be a solution of (3) with $\sigma=0$, satisfying $u=0$ on $\partial\Omega$.
Then $-\log(-u)$ is a convex function.
\end{lemma}
\emph{Proof.} Observe $\varphi:=-\log(-u)$ satisfies
$$|D\varphi|^{\frac{1}{\alpha}-1}\displaystyle{\sum_{i,j=1}^{2}}(\delta_{ij}-\frac{\varphi_{i}\varphi_{j}}{|D\varphi|^{2}})\varphi_{ij}=e^{\frac{1}{\alpha}\varphi}.$$
Since $\varphi(x)\rightarrow +\infty$ as $x\rightarrow\partial\Omega$, the result in $\cite{K}$(Theorem 3.13) implies $\varphi$ is convex. One may notice that two of the conditions required in [8] are the strict convexity of domain and the $C^{2}$ smoothness of solution. The first one can be resolved by using strictly convex domains to approximate the convex domain. For the smoothness condition, one may worry about the minimum point where the gradient vanishes and the equation is singular. Moreover, in view of the solution $u=\frac{1}{1+\alpha}|x|^{1+\alpha}$, we see when $\alpha<1$ it is not $C^{2}$ at the origin. However, by examining the proof in [8], one can see that the argument is made away from the minimum point, which means it can still be applied to our situation.

With the above lemma and the Lemma 4.4 in \cite{Wang}, we know that any convex compact ancient solution to the generalized curve shortening flow can be represented as a convex solution $u$ to equation (3) with $\sigma=0$, and if the solution to the flow sweeps the whole space, the corresponding
$u$ will be an entire solution. Thus Theorem 2 implies Corollary 1 immediately.

\emph{Remark 6.} We can also use the method in the section 4 of \cite{Wang} to construct a non-rotationally symmetric convex compact ancient solution for generalized curve shortening flow with power $\alpha\in (\frac{1}{2}, 1)$, and in fact the solution will be defined in a strip region. All we need
to do is replace Lemma 4.2, 4.3 and 4.4 in \cite{Wang} for mean curvature flow by the corresponding lemmas for the generalized curve shortening flow.

\section{\textbf{Proof of the second part of Theorem 2}}
First of all, we would like to point out that
instead of using Gage and Hamilton's exponential convergence of the curve shortening flow in \cite{GH} we need to use the corresponding exponential
convergence for the generalized curve shortening flow and we will state it as a lemma which is corresponding to lemma 3.2 in \cite{Wang}.

\begin{lemma}
Suppose $\{C_{t}\}$ be a convex solution to the generalized curve shortening flow with initial curve $\{C_{0}\}$ uniformly convex.
Suppose $\{C_{t}\}\subset N_{\delta_{0}}S^{1}$  for some unit circle $S^{1}$ , $\{C_{t}\}$ shrinks to the origin at $t=\frac{1}{1+\alpha}$. Denote $\widetilde{C_{t}}=(1-(1+\alpha)t)^{-\frac{1}{1+\alpha}}C_{t}$ as a normalization of $C_{t}$. Then
$$\widetilde{C_{t}}\subset N_{\delta_{t}}S^{1},$$
with
$$\delta_{t}\leq C\delta_{0}(\frac{1}{1+\alpha}-t)^{\iota}$$
for some small positive constant $\iota$.
\end{lemma}

The proof of the above lemma is similar to the proof of lemma 3.2 in \cite{Wang}. Using the condition that the initial curve is uniformly convex and the estimates in section II of [1], we can apply Schauder's estimates safely for $\alpha>\frac{1}{2}$ as in \cite{Wang}, which says that for $t\in (\frac{1}{4\alpha+4},\frac{1}{2\alpha+2})$, $$\|\widetilde{\ell_{t}}-S^{1}\|_{C^{k}}\leq C\delta_{0}.$$
Although the constant $C$ will depend on the lower and upper bound of the curvature of the initial curve, it is not a problem for our purpose, since when we blow down the solution for $\sigma=0$, the norm of the gradient $Du_{h}$ on the curve $\{u_{h}(x)=1\}$ approaches to 1. By the equation $\kappa u_{\gamma}^{\frac{1}{\alpha}}=1$ we see that the curvature  $\kappa$ is also very close to 1 on that curve. However, the estimates in section II of [1] also shows that when $\alpha\leq 1$ the uniformly convex condition (though the convexity is still needed) is not needed, and the constant $C$ in the above lemma is independent of the bound on the curvature of the initial curve. For the exponential decay rate of the derivatives of curvature,
one can imitate the proof in Hamilton and Gage [6], and our corresponding estimate will be $|\kappa'(\tau)|\leq C\delta_{0}e^{-\iota\tau}$ for some small positive number $\iota$, where $\tau=-\frac{1}{1+\alpha}\log(\frac{1}{1+\alpha}-t)$. This estimate immediately implies our lemma.

An alternative way to see that is by writing down the normalized evolution equation for the generalized curve shortening flow by using support function $s(\theta,\tau)$ as following $$s_{\tau}=-(s_{\theta\theta})^{-\alpha}+s,$$ here we still take the origin as the limiting point of the original generalized curve shortening flow. Then the linearized equation of the flow about the circle solution is $$s_{\tau}=\alpha(s_{\theta\theta}+s)+s.$$ The rate of convergence is governed by the eigenvalues of the right hand side.
The constant eigenfunction corresponds to scaling, which is factored out, while the $\sin\theta$ and $\cos\theta$ correspond to translations, which are also factored out.  The next is $\cos(2\theta)$, which gives eigenvalue $1-3\alpha$. So when $\alpha>\frac{1}{3},$ we have exponential convergence of the normalized solution to the limiting circle with exponent $1-3\alpha.$ The author learned this from professor Ben Andrews.

In the following we will consider the case when $\sigma=1$ .
Without loss of generality we can assume $u(0)=\inf u$.
Let $u_{h}(x)=\frac{1}{h}u(h^{\frac{1}{1+\alpha}}x)$. Then $u_{h}$ satisfies the equation $L_{\sigma}u_{h}=1$ with $\sigma=h^{-\frac{2\alpha}{1+\alpha}}.$
By Theorem 1, $\Gamma_{\frac{1}{1+\alpha},u_{h}}$ converges to the unit
circle as $h\rightarrow\infty$.

\begin{lemma}
\begin{equation}
u(x)=\frac{1}{1+\alpha}|x|^{1+\alpha}+O(|x|^{1+\alpha-2\alpha\beta})
\end{equation} where $C$ is a fixed constant and the constant $\beta$ is chosen such that $\frac{1}{2\alpha}<\beta<\min\{1,\frac{1+\alpha}{2\alpha}\}$.
\end{lemma}
For small $\delta_{0}>0$, taking $h$ large enough so that
\begin{eqnarray}
\Gamma_{\frac{1}{1+\alpha},u_{h}}\subset N_{\delta_{0}}(S^{1})
\end{eqnarray}
for unit circle $S^{1}$ with center $p_{0}$. Note that when $h$ is large, $\delta_{0}$ is very close to 0.
 Then we will prove the following claim,

\emph{ Claim 3}. For small fixed $\tau$,
\begin{eqnarray}
\Gamma_{\tau,u_{h}}\subset ((1+\alpha)\tau)^{\frac{1}{1+\alpha}}N_{\delta_{\tau}}((1+\frac{a_{0}}{\tau})^{\frac{1}{1+\alpha}}S^{1})
\end{eqnarray}
with
\begin{eqnarray}
\delta_{\tau}\leq C_{1}(\tau)\sigma^{\beta}+C_{2}\delta_{0}\tau^{\eta},
\end{eqnarray}
where the constants $C_{1}$ and $C_{2}$ are independent of $\delta_{0}$ and $h$, and $C_{2}$ is also independent of $\tau$, $\eta$ is a small
positive constant. $u_{0}$ is the solution of $L_{0}(u)=1$ in $\Omega_{\frac{1}{1+\alpha},u_{h}}$
satisfying $u_{0}=u_{h}=\frac{1}{1+\alpha}$ on $\partial\Omega_{\frac{1}{1+\alpha},u_{h}}$, $a_{0}=|\inf u_{0}|$ and the center of $(1+\frac{a_{0}}{\tau})^{\frac{1}{1+\alpha}}S^{1}$ is the minimum point of $u_{0}$ times a factor $((1+\alpha)\tau)^{-\frac{1}{1+\alpha}}$.

\emph{Proof of Claim 3.} It is equivalent to prove
\begin{eqnarray}
 \text{dist}\left((1+\alpha)^{\frac{1}{1+\alpha}}(\tau+a_{0})^{\frac{1}{1+\alpha}}S^{1},\Gamma_{\tau,u}\right)\leq C_{1}(\tau)\sigma^{\beta}+C_{2}\delta_{0}\tau^{\frac{1}{1+\alpha}+\eta},
\end{eqnarray}
where $\eta$ is some small positive constant, $C_{2}$ is independent of $\tau$.
 by Theorem 1 we know  $u_{h}$ converges to $\frac{1}{1+\alpha}|x|^{1+\alpha}$ uniformly on any compact subset of $\textbf{R}^{2}$, then by the convexity of $u_{h}$, we have that when $$x\in\{x\in\Omega_{\frac{1}{1+\alpha},u_{h}}:\tau_{0}\leq u_{h}<\frac{1}{1+\alpha}\},$$ $|Du_{h}|$
is bounded above and below by some constants depending on $\tau_{0}$ for large $h$, by the growth condition for $D^{2}u$ in Theorem 2 we have $\sigma(u_{h})_{\gamma\gamma}\leq C\sigma^{\beta},$ where $C$ is a constant depending on $\tau_{0}$. Therefore we have $\kappa (u_{h})_{\gamma}^{\frac{1}{\alpha}}\approx 1-C\sigma^{\beta}\ \text{on}\ \{x\in\Omega_{\frac{1}{1+\alpha},u_{h}}:\tau\leq u_{h}<\frac{1}{1+\alpha}\},$
 where $C$ depends on $\tau_{0}$. Denote $$\widetilde{u}_{0}=(1-C\sigma^{\beta})^{\alpha}(u_{0}-\frac{1}{1+\alpha})+\frac{1}{1+\alpha},$$ then $$L_{0}(\widetilde{u}_{0})=1-C\sigma^{\beta}\ \text{in}\ \Omega_{\frac{1}{1+\alpha},u_{h}}$$
with $\widetilde{u}_{0}=u_{h}=\frac{1}{1+\alpha}$ on $\partial\Omega_{\frac{1}{1+\alpha},u_{h}}$.
Now by comparison principle we have
$\Omega_{\tau,u_{0}}\subset\Omega_{\tau,u_{h}}\subset\Omega_{\tau,\widetilde{u}_{0}}$, and by the asymptotic behavior of $u_{0}$ we have
$$\Gamma_{\tau,u_{0}}\subset N_{\zeta}((\tau+a_{0})^{\frac{1}{1+\alpha}}S^{1})\ \text{and}\ \Gamma_{\tau, \widetilde{u}_{0}}\subset N_{\zeta}((\tau+a_{0}-C\sigma^{\beta})^{\frac{1}{1+\alpha}}S^{1}),$$ where $\zeta={C\delta_{0}(\tau+a_{0})^{\eta}}$. Denote $\ell_{1}=(\tau+a_{0})^{\frac{1}{1+\alpha}}S^{1}$, $\ell_{2}=(\tau+a_{0}-C\sigma^{\beta})^{\frac{1}{1+\alpha}}S^{1}$, both of them are centered at $p_{1}$,
which is the minimum point of $u_{0}$. Hence $\text{dist}((\tau+a_{0})^{\frac{1}{1+\alpha}}S^{1},\Gamma_{\tau,u_{h}})\leq\text{dist}(\ell_{1},\ell_{2})+C\delta_{0}(\tau+a_{0})^{\frac{1}{1+\alpha}+\eta}$, where $\text{dist}(\ell_{1},\ell_{2})$ can be bounded by $C_{1}(\tau)\sigma^{\beta}$, hence (26) follows from
the above discussion. Now we will use an iteration argument to prove the following Claim 4, which will enable us to simplify (25) and (26).

\emph{Claim 4}:
\begin{equation}
a_{0}\leq\begin{cases}
C\sigma|\log(\sigma)| & \text{if}\ \alpha\leq 1 \\
C\sigma^{\frac{1+\alpha}{2\alpha}} & \text{if}\ \alpha>1
\end{cases}
\end{equation}
\emph{Proof of Claim 4.} We fix a large constant $A$ such that $\{u_{\frac{A}{\tau}}=\frac{1}{1+\alpha}\}$ is very close to a unit circle. Let
$u_{0,\tau^{k}}$ solve $L_{0}u=1$ with boundary condition $u=\tau^{k}$ on $\{u_{h}=\tau^{k}\}$. Denote $a_{k}=|\inf u_{0,\tau^{k}}|.$ From the proof of Claim 3 we see that $\{u_{0}<\tau\}\supset \{u_{0,\tau}<\tau\}\supset \{\widetilde{u}_{0}<\tau\},$ by comparison principle, we have
$\inf u_{0}<\inf u_{0,\tau}<\inf \widetilde{u}_{0}.$ So by the construction of $\widetilde{u}_{0}$ and a simple computation, we have $a_{0}-a_{1}\leq \inf \widetilde{u}_{0}-\inf u_{0}\leq C\sigma.$ When $\tau^{k}\geq \frac{A}{h}$, we can iterate this argument for $u_{0,\tau^{k}}$ and $u_{0,\tau^{k+1}}$ by rescaling them to $\frac{1}{1+\alpha}\tau^{-k}u_{0,\tau^{k}}\left((1+\alpha)^{\frac{1}{1+\alpha}}\tau^{\frac{k}{1+\alpha}}x\right)$ and
$\frac{1}{1+\alpha}\tau^{-k}u_{0,\tau^{k+1}}\left((1+\alpha)^{\frac{1}{1+\alpha}}\tau^{\frac{k}{1+\alpha}}x\right)$ respectively, after rescaling back, we have $a_{k}-a_{k+1}\leq C\sigma$. Note that the choice of $A$ and the condition $\tau^{k}\geq \frac{A}{h}$ ensure the uniform gradient bound needed in the above argument. Let $k_{0}$ be an integer satisfying $\tau^{k_{0}}\geq \frac{A}{h}\geq \tau^{k_{0}+1}$, after $k_{0}$ steps we stop the iteration, and notice that $\{u_{h}=\frac{A}{h}\}=\frac{1}{h^{\frac{1}{1+\alpha}}}\{u=A\}$ is contained in a circle with radius $Ch^{-\frac{1}{1+\alpha}}$ for some constant $C$, so it takes at most time $Ch^{-1}=C\sigma^{\frac{1+\alpha}{2\alpha}}$ for $\{u_{h}=\frac{A}{h}\}$ shrink into a point. Claim 4 follows from the above discussion.

By omitting the lower order term we can rewrite (25) and (26) as
\begin{eqnarray*}
\Gamma_{\tau,u_{h}}\subset ((1+\alpha)\tau)^{\frac{1}{1+\alpha}}N_{\delta_{\tau}}(S^{1})
\end{eqnarray*}
with
\begin{equation}
\delta_{\tau}\leq C_{1}(\tau)\sigma^{\beta}+C_{2}\delta_{0}\tau^{\eta}.
\end{equation}

If we take $\tau$ small such that $C_{2}\tau^{\eta}\leq \frac{1}{4}$, (29) becomes
\begin{equation}
\delta_{\tau}\leq
C_{1}(\tau)\sigma^{\beta}+\frac{1}{4}\delta_{0}.
\end{equation}
Now we can carry out an iteration argument similar as that in \cite{Wang}. We start at the level $\frac{1}{1+\alpha}\tau^{-k_{0}}$ for $k_{0}$ very large.
Denote $\Omega_{k}=\tau^{\frac{k}{1+\alpha}}\Omega_{\frac{1}{1+\alpha}\tau^{-k},u}$ and $\Gamma_{k}=\partial \Omega_{k}$. Define $\delta_{k}$ similarly to that in $\cite{Wang}$.
 By (30) we have
\begin{equation}
\delta_{k-1}\leq
C_{1}(\tau)\tau^{(k-1)\frac{2\alpha\beta}{1+\alpha}}+\frac{1}{4}\delta_{k}
\end{equation}
for $k=k_{0},k_{0}+1,\cdots$.
Then we have
\begin{eqnarray}
\Gamma_{j}\subset N_{\delta_{j}}(S^{1})
\end{eqnarray}
with
\begin{equation}
\delta_{j}\leq
C\tau^{j\frac{2\alpha\beta}{1+\alpha}}
\end{equation}
 It follows that
\begin{eqnarray}
\Gamma_{\frac{1}{1+\alpha}\tau^{-j},u}\subset N_{\widetilde{\delta}_{j}}(\tau^{\frac{-j}{1+\alpha}}S^{1})
\end{eqnarray}
with
\begin{equation}
\widetilde{\delta}_{j}\leq
C\tau^{\frac{2\alpha\beta-1}{1+\alpha}j}
\end{equation}
where $\tau^{\frac{-j}{1+\alpha}}S^{1}$
is centered at $z_{j}=\tau^{\frac{-j}{1+\alpha}}y_{j}$.
From Lemma 3 and (30) it is not hard to see that we have
\begin{equation}
|z_{j}-z_{j-1}|\leq C\tau^{\frac{2\alpha\beta-1}{1+\alpha}j}
\end{equation}
Denote $z_{0}=\lim_{j\rightarrow\infty}z_{j}$.
Then \begin{equation}
|z_{j}-z_{0}|\leq C\tau^{\frac{2\alpha\beta-1}{1+\alpha}j},
\end{equation} which means in (34) we can assume the circle is centered at $z_{0}$ by changing the constant $C$ a little bit. In fact when we choose different $\tau$, the corresponding $z_{0}$ will not change, so we can assume $z_{0}=0$. Hence for
$h=\frac{1}{1+\alpha}\tau^{-j}$,
$$\Gamma_{h,u}\subset N_{\delta}\left((1+\alpha)^{\frac{1}{1+\alpha}}h^{\frac{1}{1+\alpha}}S^{1}\right),$$
where
\begin{equation}
\delta\leq
Ch^{\frac{1-2\alpha\beta}{1+\alpha}}
\end{equation}
and $S^{1}$ is centered at the origin. By using different $\tau$, it is easy to see that estimate holds for any large $h$. Lemma 5 follows from the above estimates.

Now we can finish the proof of Theorem 2 in the following way.

\emph{Proof of the second part of Theorem 2.} Denote $u^{*}$ as the Legendre transform of $u$. Then $u^{*}$ satisfies the following equation
\begin{eqnarray}
G(x,D^{2}u^{*})=\frac{\det{D^{2}u^{*}}}{(\delta_{ij}-\frac{x_{i}x_{j}}{1+|x|^{2}})F^{ij}(u^{*})}=(1+|x|^{2})^{\frac{1}{2\alpha}-\frac{1}{2}},
\end{eqnarray}
where $F^{ij}(u^{*})=\frac{\partial \det{r}}{\partial r_{ij}}$, at  $r=D^{2}(u^{*}).$
We have \begin{eqnarray}
u^{*}(x)=C(\alpha)|x|^{1+\alpha}+O(|x|^{\frac{1+\alpha-2\alpha\beta}{\alpha}}),
\end{eqnarray}
where $C(\alpha)$ is a constant depending only
on $\alpha$.
In fact, for big $h$, by Lemma 5 we have $$u_{h}(x)=\frac{1}{1+\alpha}|x|^{1+\alpha}+O(|h|^{\frac{-2\alpha\beta}{1+\alpha}})$$ in $B_{1}(0).$
Denote $u_{h}^{*}$ as the Legendre transforms of $u_{h}$. Then $$u_{h}^{*}(x)=C(\alpha)|x|^{1+\frac{1}{\alpha}}+O(|h|^{\frac{-2\alpha\beta}{1+\alpha}}),$$ where $C(\alpha)$ is a constant depending only
on $\alpha$ and in fact it is comes from the Legendre transform of the function $\frac{1}{1+\alpha}|x|^{1+\alpha}.$ Note that
$u_{h}^{*}(x)=h^{-1}u^{*}(h^{\frac{\alpha}{1+\alpha}}x)$, we obtain (40).

Let $u_{0}$ be the unique radial solution of (3) with $\sigma=1$, and let $u_{0}^{*}$ be the Legendre transform of $u_{0}$. Similar to (40)
we have
\begin{eqnarray}
u_{0}^{*}(x)=C(\alpha)|x|^{1+\alpha}+O(|x|^{\frac{1+\alpha-2\alpha\beta}{\alpha}}).
\end{eqnarray}
Since both $u^{*}$ and $u_{0}^{*}$ satisfy equation (39), $v=u^{*}-u_{0}^{*}$ satisfies the following elliptic equation
$$\displaystyle{\sum_{i,j=1}^{n}}a_{ij}(x)v_{ij}=0\ \text{in}\ \textbf{R}^{2},$$ where
$$a_{ij}=\int_{0}^{1}G^{ij}(x, D^{2}u_{0}^{*}+t(D^{2}u^{*}-D^{2}u_{0}^{*})dt,$$ here $G^{ij}=\frac{\partial G(x,r)}{\partial r_{ij}}$ for any
symmetric matrix $r$. Note that by the choice of $\beta$, $\frac{1+\alpha-2\alpha\beta}{\alpha}<1$, so by (40) and (41) $v=O(|x|^{\frac{1+\alpha-2\alpha\beta}{\alpha}})=o(|x|)$, as $|x|\rightarrow\infty$. Using the Liouville Theorem by Bernstein \cite{Simon} (p.245)
we conclude that $v$ is a constant.

\end{document}